\documentclass{amsart}
\usepackage[utf8]{inputenc}
\usepackage{cases, url}
\usepackage[all]{xy}

\newtheorem{Thm}{Theorem}[section]
\newtheorem{Lem}[Thm]{Lemma}
\newtheorem{Cor}[Thm]{Corollary}
\newtheorem{Prop}[Thm]{Proposition}

\theoremstyle{definition}
\newtheorem{Rem}[Thm]{Remark}

\newcommand{\CC}{\mathbb{C}}

\newcommand{\ZZ}{\mathbb{Z}}

\newcommand{\NN}{\mathbb{N}}

\newcommand{\cB}{{\mathcal B}}

\newcommand{\cF}{{\mathcal F}}

\newcommand{\cI}{{\mathcal I}}

\newcommand{\cM}{{\mathcal M}}

\newcommand{\cO}{{\mathcal O}}
\newcommand{\cP}{{\mathcal P}}

\newcommand{\cS}{{\mathcal S}}

\newcommand{\fg}{\mathfrak{g}}

\newcommand{\fn}{\mathfrak{n}}

\newcommand{\bd}{\mathbf{d}}
\newcommand{\be}{\mathbf{e}}
\newcommand{\bg}{\mathbf{g}}
\newcommand{\bi}{\mathbf{i}}
\newcommand{\bj}{\mathbf{j}}
\newcommand{\bm}{\mathbf{m}}
\newcommand{\bn}{\mathbf{n}}
\newcommand{\bp}{\mathbf{p}}
\newcommand{\bq}{\mathbf{q}}
\newcommand{\br}{\mathbf{r}}
\newcommand{\bE}{\mathbf{E}}
\newcommand{\bM}{\mathbf{M}}
\newcommand{\bU}{\mathbf{U}}
\newcommand{\bX}{\mathbf{X}}

\newcommand{\Del}{\Delta}
\newcommand{\Gam}{\Gamma}

\newcommand{\Ome}{\Omega}
\newcommand{\alp}{\alpha}
\newcommand{\bet}{\beta}
\newcommand{\gam}{\gamma}
\newcommand{\del}{\delta}   
\newcommand{\eps}{\epsilon}  
\newcommand{\lam}{\lambda}

\newcommand{\vph}{\varphi}
\newcommand{\vth}{\theta}

\makeatletter
\let\save@mathaccent\mathaccent
\newcommand*\if@single[3]{%
  \setbox0\hbox{${\mathaccent"0362{#1}}^H$}%
  \setbox2\hbox{${\mathaccent"0362{\kern0pt#1}}^H$}%
  \ifdim\ht0=\ht2 #3\else #2\fi
  }
\newcommand*\rel@kern[1]{\kern#1\dimexpr\macc@kerna}
\newcommand*\widebar[1]{\@ifnextchar^{{\wide@bar{#1}{0}}}{\wide@bar{#1}{1}}}
\newcommand*\wide@bar[2]{\if@single{#1}{\wide@bar@{#1}{#2}{1}}{\wide@bar@{#1}{#2}{2}}}
\newcommand*\wide@bar@[3]{%
  \begingroup
  \def\mathaccent##1##2{%
    \let\mathaccent\save@mathaccent
    \if#32 \let\macc@nucleus\first@char \fi
    \setbox\z@\hbox{$\macc@style{\macc@nucleus}_{}$}%
    \setbox\tw@\hbox{$\macc@style{\macc@nucleus}{}_{}$}%
    \dimen@\wd\tw@
    \advance\dimen@-\wd\z@
    \divide\dimen@ 3
    \@tempdima\wd\tw@
    \advance\@tempdima-\scriptspace
    \divide\@tempdima 10
    \advance\dimen@-\@tempdima
    \ifdim\dimen@>\z@ \dimen@0pt\fi
    \rel@kern{0.6}\kern-\dimen@
    \if#31
      \overline{\rel@kern{-0.6}\kern\dimen@\macc@nucleus\rel@kern{0.4}\kern\dimen@}%
      \advance\dimen@0.4\dimexpr\macc@kerna
      \let\final@kern#2%
      \ifdim\dimen@<\z@ \let\final@kern1\fi
      \if\final@kern1 \kern-\dimen@\fi
    \else
      \overline{\rel@kern{-0.6}\kern\dimen@#1}%
    \fi
  }%
  \macc@depth\@ne
  \let\math@bgroup\@empty \let\math@egroup\macc@set@skewchar
  \mathsurround\z@ \frozen@everymath{\mathgroup\macc@group\relax}%
  \macc@set@skewchar\relax
  \let\mathaccentV\macc@nested@a
  \if#31
    \macc@nested@a\relax111{#1}%
  \else
    \def\gobble@till@marker##1\endmarker{}%
    \futurelet\first@char\gobble@till@marker#1\endmarker
    \ifcat\noexpand\first@char A\else
      \def\first@char{}%
    \fi
    \macc@nested@a\relax111{\first@char}%
  \fi
  \endgroup
}
\makeatother

\newcommand{\bQ}{\widebar{Q}}

\newcommand{\oH}{\widebar{H}}
\newcommand{\oB}{\widebar{B}}
\newcommand{\oOme}{\widebar{\Ome}}
\newcommand{\bcM}{\widebar{\cM}}

\newcommand{\ud}{\underline{d}}

\newcommand{\tth}{\tilde{\theta}}

\newcommand{\id}{\operatorname{id}}
\newcommand{\ad}{\operatorname{ad}}
\newcommand{\diag}{\operatorname{diag}}
\newcommand{\sgn}{\operatorname{sgn}}
\newcommand{\rk}{\operatorname{rank}}
\newcommand{\dimv}{\underline{\dim}}
\newcommand{\rkv}{\underline{\rk}}
\newcommand{\supp}{\operatorname{supp}}
\newcommand{\lcm}{\operatorname{lcm}}

\newcommand{\rep}{\operatorname{rep}} 
\newcommand{\Rep}{\operatorname{Rep}} 
\newcommand{\Gr}{\operatorname{Gr}}
\newcommand{\Grlf}{\operatorname{Grlf}}
\newcommand{\Fl}{\operatorname{Fl}}

\newcommand{\pdim}{\operatorname{proj.dim}}
\newcommand{\idim}{\operatorname{inj.dim}}

\newcommand{\sub}{\operatorname{sub}}
\newcommand{\fac}{\operatorname{fac}}

\newcommand{\Hom}{\operatorname{Hom}}
\newcommand{\Ext}{\operatorname{Ext}}
\newcommand{\End}{\operatorname{End}}
\newcommand{\Aut}{\operatorname{Aut}}
\newcommand{\GL}{\operatorname{GL}}

\newcommand{\Ima}{\operatorname{Im}}
\newcommand{\Ker}{\operatorname{Ker}}

\newcommand{\Irr}{\operatorname{Irr}}

\newcommand{\crys}{\mathrm{cr}}
\newcommand{\op}{\mathrm{op}}
\newcommand{\out}{\mathrm{out}}
\newcommand{\inn}{\mathrm{in}}
\newcommand{\real}{\mathrm{re}}
\newcommand{\imag}{\mathrm{im}}
\newcommand{\lfr}{\mathrm{l.f.}}
\newcommand{\fib}{\mathrm{fib}}

\newcommand{\cyc}{\mathrm{cyc}}

\newcommand{\wt}{\operatorname{wt}}

\newcommand{\df}{\colon}

\newcommand{\Iff}{\Longleftrightarrow}
\newcommand{\ra}{\rightarrow}
\newcommand{\bil}[1]{\langle #1\rangle}
\newcommand{\abs}[1]{\left| #1\right|}

\title[Quivers with relations for symmetrizable Cartan matrices]%
{Quivers with relations for symmetrizable Cartan matrices and algebraic
Lie theory}

\author{Christof Geiß}
\address{Instituto de Matemáticas,
Universidad Nacional Autónoma de México, 
Ciudad Universitaria,
04510 Cd. de México, MEXICO}
\email{christof.geiss@im.unam.mx}


\begin{document}
\begin{abstract}
We give an overview of our effort to introduce (dual) semicanonical bases in the
setting of symmetrizable Cartan matrices.
\end{abstract}

\numberwithin{equation}{section}

\maketitle

\section{Introduction}
One of the original motivations of Fomin and Zelevinsky for introducing
cluster algebras was ``to understand, in a concrete and combinatorial way, G.
Lusztig's theory of total positivity and canonical bases''~\cite{Fo}. 
This raised the question of finding a cluster algebra structure on the 
coordinate ring
of a unipotent cell, and to study its relation with Lusztig's bases. In a
series of works culminating with~\cite{gls1} and~\cite{gls2}, we showed that the
coordinate ring of a unipotent cell of a symmetric Kac-Moody group has
indeed a cluster algebra structure, whose cluster monomials belong to the
dual of Lusztig's semicanonical basis of the enveloping algebra of the
attached Kac-Moody algebra.  
Since the semicanonical basis is built in
terms of constructible functions on the complex varieties of nilpotent
representations of the preprojective algebra of a quiver, it is not 
straightforward  to extend those results to the  setting of symmetrizable Cartan
matrices, which appears more natural from the Lie theoretic point of view.
The purpose of these notes is to give an overview of~\cite{qrel1} - 
\cite{qrel5}, where we are trying to make progress into this direction.

The starting point of our project was~\cite{HL}, where
Hernandez and Leclerc observed that certain quivers with potential allowed
to encode the q-characters of the Kirillov-Reshetikhin modules of the quantum 
loop algebra $U_q(L\fg)$, where $\fg$ is a complex simple Lie algebra of
arbitrary Dynkin type. This quiver with potential served as model for the
definition of our generalized preprojective algebras $\Pi=\Pi_K(C,D)$ associated
to a symmetrizable Cartan matrix $C$ with symmetrizer $D$ over an arbitrary 
field $K$, which extends the classical construction of Gelfand and 
Ponomarev~\cite{GP}.
After the completion of a preliminary version of~\cite{qrel1} we 
learned that Cecotti and Del~Zotto~\cite{CDZ} and Yamakawa~\cite{Yam} had
introduced similar constructions for quite different reasons. 
In comparison to the classical constructions
of Dlab and Ringel~\cite{DR1}, \cite{DR3} for a symmetrizable 
Cartan matrix $C$, we replace field
extensions by truncated polynomial rings. Many of the core results of
representations of species carry over over to this setting if we restrict
our attention to the so-called locally free modules, see~\cite{qrel1}.
In particular, we have for each orientation $\Ome$ of $C$ an algebra 
$H=H_K(C,D,\Ome)$ such that in many respects $\Pi$ can be considered as the
preprojective algebra of $H$. Our presentation of these results 
in Section~\ref{sec:qrel1} is inspired by the thesis~\cite{Geu}, where 
Geuenich obtains similar results for a larger class of algebras.

Since our construction works in particular over algebraically closed fields,  
we can extend to our algebras $H$ and $\Pi$ several basic results about 
representation varieties of  quivers and of varieties of nilpotent 
representations of the preprojective  algebra of a  quiver in our new context, 
again if we restrict our attention to locally free modules, see 
Section~\ref{sec:rvar}.
Nandakumar and Tingley~\cite{NT} obtained similar results by
studying the set of $K$-rational points of the representation scheme of a
species preprojective algebra, which is defined over certain infinite,
non algebraically closed  fields $K$.

In our setting we can take $K=\CC$, and study algebras of constructible 
functions on those varieties of locally free modules and 
realize in this manner
the universal enveloping algebra $U(\fn)$ of the positive part $\fn$ of
a complex semisimple Lie algebra, together with a Ringel type PBW-basis in
terms of the representations of $H$. 
For arbitrary symmetrizable Cartan matrices we
can realize $U(\fn)$ together with a semicanonical basis, modulo our
\emph{support conjecture}, see Section~\ref{sec:const}.\\[1ex]
\textbf{Conventions.} We use basic concepts from representation
theory of finite dimensional algebras, like Auslander-Reiten theory or
tilting theory without further reference. A good source for this material
is~\cite{Ri3}. For us, a quiver is an oriented graph $Q=(Q_0,Q_1,s,t)$ with
vertex set $Q_0$, arrow set $Q_1$ and functions $s,t\df Q_1\ra Q_0$ indicating
the start and terminal point of each arrow. We also write $D=\Hom_K(-,K)$.
We say that an $A$-module $M$ is \emph{rigid} if $\Ext^1_A(M,M)=0$.

\section{Combinatorics of symmetrizable Cartan matrices}
\label{sec:comb}
\subsection{Symmetrizable Cartan matrices and quivers}
Let $I=\{1,2,\ldots,n\}$.  
A \emph{symmetrizable  Cartan matrix} is 
an integer matrix  $C=(c_{ij})\in\ZZ^{I\times I}$ such that the following holds: 
\begin{itemize}
\item $c_{ii}=2$ for all $i\in I$ and 
 $c_{ij}\leq 0$ for all $i\neq j$,
\item there exist $(c_i)_{i\in I}\in\NN_+^I$ such that 
$\diag(c_1,\ldots,c_n)\cdot C$ is a symmetric.
\end{itemize}
In this situation $D:=\diag(c_1,\ldots,c_n)\in\ZZ^{I\times I}$ is called the 
\emph{symmetrizer} of $C$. 
Note that the symmetrizer is not unique. In particular, for all
$k\in\NN_+$ also $kD$ is a symmetrizer of $C$.

It is easy to see that the datum $(C,D)$ of a symmetrizable Cartan matrix
$C$ and its symmetrizer $D$ 
is equivalent to displaying
a weighted graph $(\Gam,\ud)$ with 
\begin{itemize}
\item $I$ the set of vertices of $\Gam$,
\item $g_{ij}:=\gcd(c_{ij},c_{ji})$ edges between $i$ and $j$,
\item $\ud\df I\ra\NN_+, i\mapsto c_i$.
\end{itemize}
Here we agree that $\gcd(0,0)=0$. We have then 
$c_{ij}=-\frac{\lcm(c_i,c_j)}{c_i} g_{ij}$ for all $i\neq j$.

\subsection{Bilinear forms,  reflections and roots}\label{ssec:root}
We identify the root lattice of the Kac-Moody Lie algebra $\fg(C)$
associated to $C$ with $\ZZ^I =\oplus_{i\in I} \ZZ\alp_i$,
where the simple roots $(\alpha_i)_{i\in I}$ form the standard basis. 
We define on $\ZZ^I$ by 
\[
(\alp_i, \alp_j)_{C,D}=c_ic_{ij},
\]
a symmetric bilinear form. 
The \emph{Weyl group} $W=W(C)$
is the subgroup of $\Aut(\ZZ^I)$, which is generated
by the simple reflections $s_i$ for $i\in I$, where
\[
s_i(\alp_j) = \alp_j-c_{ij}\alp_i.
\]
The \emph{real roots} are the set
\[
\Del_\real(C):=\cup_{i\in I} W(\alp_i).
\]
The \emph{fundamental region} is
\[
F:=\{\alp\in\NN^I\mid\supp(\alp)\text{ is connected, and } 
(\alp,\alp_i)_{C,D}\leq 0
\text{ for all } i\in I\}.
\]
Here, $\supp(\alp)$ is the full subgraph of $\Gam(C)$ with vertex
set $\{i\in I\mid \alp(i)\neq 0\}$. 
Then the imaginary roots are by definition the set
\[
\Del_\imag(C):= W(F)\cup W(-F).
\]
Finally the set of all roots is
\[
\Del(C):=\Del_\real\cup\Del_\imag(C).
\]
The positive roots are
$\Del^+(C) :=\Del(C)\cap\NN^I$, and 
it is remarkable that
$\Del(C)= \Del^+(C)\cup -\Del^+(C)$.

A sequence $\bi=(i_1, i_2,\ldots, i_l)\in I^l$ is called a 
\emph{reduced expression} for $w\in W$ if $w=s_{i_l}\cdots s_{i_2} s_{i_1}$
and $w$ can't be expressed as a product of less than $l=l(w)$ 
reflections of the form $s_i$ $(i\in I)$. In this case we set
\begin{equation} \label{eq:rootlist}
\bet_{\bi,k}:=s_{i_1}s_{i_2}\cdots s_{i_{k-1}}(\alp_{i_k}) \text{ and }
\gam_{\bi,k}:=s_{i_l}s_{i_{l-1}}\cdots s_{i_{k+1}}(\alp_{i_k})
\end{equation}
for $k=1,2,\ldots, l$, and understand $\bet_{\bi,1}=\alp_{i_1}$ as well as 
$\gam_{\bi,l}=\alp_{i_l}$. It is a standard fact that 
$\bet_{\bi,k}\in\Del^+$ for $k=1,2,\ldots, l$, and that these roots are pairwise 
different. Obviously,
\[
w(\bet_{\bi,k})=-\gam_{\bi,k}\text{ for } k=1,2,\ldots, l.
\]
The following result is well known.
\begin{Prop} \label{prp:Dynkin}
For a connected, symmetrizable Cartan matrix $C$ the following
are equivalent:
\begin{itemize}
\item $C$ is of Dynkin type.
\item The Weyl group $W(C)$ is finite.
\item The root system $\Del(C)$ is finite
\item All roots are real: $\Del(C)=\Del_\real(C)$.
\end{itemize}
Moreover, if in this situation $\bi$ is a reduced expression for $w_0$, 
the longest element
of $W$, then $\Del^+=\{\bet_{\bi,1},\bet_{\bi,2},\ldots,\bet_{\bi,l}\}$.
\end{Prop}

\subsection{Orientation and Coxeter elements}\label{ssec:cox}
An \emph{orientation} of $C$ is a set $\Ome\subset I\times I$ such that
\begin{itemize}
\item $\abs{\Ome\cap\{(i,j),(j,i)\}}\Iff c_{ij}<0$,
\item for each sequence $i_1, i_2,\ldots,i_{k+1}$  with 
$(i_j,i_{j+1})\in\Ome$ for  $j=1,2,\ldots,k$ we have $i_1\neq i_{k+1}$.
\end{itemize}
The orientation $\Ome$ can be interpreted as upgrading 
the weighted graph $(\Gam,\ud)$ of $(C,D)$ to a weighted quiver $(Q^\circ,\ud)$  
with $g_{ij}$ arrows 
$\alp_{ij}^{(1)},\ldots \alp_{ij}^{(g_{ij})}$ from $j$ to $i$ if $(i,j)\in\Ome$,
such that  $Q^{\circ}=Q^{\circ}(C,\Ome)$ has no oriented cycles.

For an orientation $\Ome$ of the symmetrizable Cartan matrix 
$C\in\ZZ^{I\times I}$ and $i\in I$ we define
\[
s_i(\Ome):=\{(r,s)\in\Ome\mid i\not\in\{r,s\}\}\cap
\{(s,r)\in I\times I\mid (r,s)\in\Ome\text{ and } i\in\{r,s\}\}.
\]
Thus, in $Q^\circ(C,s_i(\Ome))$  the orientation of precisely the arrows in 
$Q^\circ(C,\Ome)$,  which are incident with $i$, is changed.
If $i$ is a sink or a source of $Q^\circ(C,\Ome)$ then $s_i(\Ome)$ is also
an orientation of $C$. It is convenient to define
\[
\Ome(-,i):=\{j\in I\mid (j,i)\in\Ome\} \text{ and }
\Ome(j,-):=\{i\in I\mid (j,i)\in\Ome\}.
\]
We have on $\ZZ^I$ the \emph{non-symmetric bilinear form}
\begin{equation} \label{eq:non-symbil}
\bil{-,-}_{C,D,\Ome}\df\ZZ^I\times\ZZ^I\ra\ZZ, 
(\alp_i,\alp_j)\mapsto\begin{cases}
c_i      &\text{if } i=j,\\
c_ic_{ij} &\text{if } (j,i)\in\Ome,\\
0        &\text{else.}
\end{cases} 
\end{equation}
We leave it as an exercise to verify that
\begin{equation} \label{eq:refl-ring}
\bil{\alp,\bet}_{C,D,\Ome}=\bil{s_i(\alp),s_i(\bet)}_{C,D, s_i(\Ome)} 
\end{equation}
if $i$ is a sink or a source for $\Ome$.

We say that a reduced expression $\bi=(i_1,i_2,\ldots,i_l)$ of $w\in W$ is
$+$-\emph{admissible} for $\Ome$ if $i_1$ is a sink of $Q^\circ(C,\Ome)$,
and $i_k$ is a sink of $Q^\circ(C,s_{i_{k-1}}\cdots s_{i_2}s_{i_1}(\Ome))$ for
$k=2,3,\ldots,l$. 
If moreover $l=n$ and $\{i_1,\ldots,i_n\}=I$, we say
that $c= s_{i_n}\cdots s_{i_2}s_{i_1}$ is the \emph{Coxeter element}
for $(C,\Ome)$. 

\subsection{Kac-Moody Lie algebras}
For a symmetrizable Cartan matrix $C\in\ZZ^{I\times I}$, the derived 
Kac-Moody Lie algebra $\fg'=\fg'(C)$ over the complex numbers has a 
presentation by $3n$ generators $e_i, h_i, f_i$ 
$(i\in I)$ subject to the following relations:
\begin{itemize}
 \item[(i)]
 $[e_i,f_j] = \del_{ij}h_i$;
 \item[(ii)]
 $[h_i,h_j]=0$; 
 \item[(iii)]
 $[h_i,e_j]=c_{ij}e_j,\quad [h_i,f_j]=-c_{ij}f_j$;
 \item[(iv)] 
 $(\ad e_i)^{1-c_{ij}}(e_j) = 0, \quad (\ad f_i)^{1-c_{ij}}(f_j) = 0$\qquad $(i\not = j)$.
\end{itemize} 
Note that for $C$ of  Dynkin type this is the Serre presentation of the
corresponding semisimple Lie algebra. In case $\rk C < \abs{I} $ we have 
of $\fg'(C)\neq\fg(C)$ and the latter has in this case a 
slightly larger Cartan subalgebra, which makes for a more complicated 
definition, see for example~\cite[Sec.~5.1]{qrel4} for a few more details. 
Of course, the main reference is~\cite{Ka}.

Let $\fn = \fn(C)$ be the Lie subalgebra generated by the $e_i\ (i\in I)$.
Then $U(\fn)$ is the associative $\CC$-algebra with generators 
$e_i\ (1\le i\le n)$ subject to the relations
\begin{equation} \label{eq:SerRel}
 (\ad e_i)^{1-c_{ij}}(e_j)  
= 0, \qquad (i,j\in I, i\neq j).
\end{equation}
$U(\fn)$ is $\NN^I$ graded with $\deg(e_i)=\alp_i$ ($i\in I$). With
\[
\fn_\alp:=\fn\cap U(\fn)_\alp\text{ for } \alp\in\Del^+(C)
\]
we recover the usual root space decomposition of $\fn$.

\section{Quivers with relations for symmetrizable Cartan matrices}
\label{sec:qrel1}
We keep the notations from the previous section, in particular 
$C\in\ZZ^{I\times I}$ is a symmetrizable Cartan matrix with symmetrizer
$D$ and $\Ome$ is an orientation for $C$.
 
\subsection{A class of 1-Iwanaga-Gorenstein algebras}
Let $K$ be a field and $Q=Q(C,D,\Ome)$ the quiver  obtained from
$Q^\circ(C,D,\Ome)$, see Section~\ref{ssec:cox}, by adding a loop $\eps_i$ 
at each vertex $i\in I$.
Then $H=H_K(C,D,\Ome)$ is the path algebra $KQ$ modulo the ideal which
is generated by the following relations:
\begin{itemize}
\item $\eps_i^{c_i}$ for all $i\in I$
\item  $\eps_i^{-c_{ji}/g_{ji}}\alp_{ij}^{(k)}-\alp_{ij}^{(k)}\eps_j^{-c_{ij}/g_{ij}}$
for all $(i,j)\in\Ome$ and $k=1,2,\ldots g_{ij}$.
\end{itemize}
Recall that $g_{ij}=g_{ji}=\gcd(c_{ij},c_{ji})$, thus 
$-c_{ij}/g_{ij}=\lcm(c_i,c_j)/c_i$.

For $(i,j)\in\Ome$ let $c'_{ij}=c_{ij}/g_{ij}$ and 
$c'_{ji}=c_{ji}/g_{ij}$. We may consider the following symmetrizable 
Cartan matrix,  symmetrizer and orientation:
\[
C^{(i,j)}=\begin{pmatrix} 2& c'_{ij}\\c'_{ji}& 2\end{pmatrix},\quad 
D^{(i,j)}=\begin{pmatrix} c_i& 0\\0 & c_j\end{pmatrix}\quad\text{and}\quad
\Ome^{(i,j)}=\{(i,j)\}.
\]
Thus, 
\[
Q^{(i,j)}:=Q(C^{(i,j)}, \Ome^{(i,j)})
\quad=\quad\xymatrix{i \ar@(ul,dl)[]|{\eps_i}&
\ar[l]_{\alp_{ij}}j\ar@(ur,dr)[]|{\eps_j}}
\]
and
\[
H^{(i,j)}:=H_K(C^{(i,j)}, D^{(i,j)},\Ome^{(i,j)})=
KQ^{(i,j)}/\bil{\eps_i^{c_i},\eps_j^{c_j},
\eps_i^{-c'_{ji}}\alp_{ij}-\alp_{ij}\eps_j^{-c'_{ij}}}.
\]
Note, that with 
\[
{_iH'_j}:= e_iH^{(i,j)} e_j\quad\text{and}\quad 
H_i=e_iH_ie_i = K[\eps_i]/(\eps_i^{c_i})
\]
it is easy to see that ${_iH_j}:={_iH'_j}^{\oplus g_{ij}}$ 
is a $H_i$-$H_j$-bimodule, 
which is free of rank $-c_{ij}$ as a $H_i$-module, and free of rank $-c_{ji}$
as $H_j$-(right)-module. 
If we define similarly $H^{(j,i)}:=H_K(C^{(i,j)},D^{(i,j)},\{(j,i)\})$ and
${_jH'_i}:=e_jH^{(j,i)}e_j$, then ${_jH_i}={_jH'_i}^{\oplus g_{ij}}$ is a 
$H_j$-$H_i$-bimodule, which is free of rank $-c_{ji}$ as $H_j$-module and free 
of rank $-c_{ij}$ as $H_i$-(left)-module. It is easy  to see that we get an
isomorphism of $H_i$-$H_j$-bimodules
\[
{_iH_j}\cong \Hom_K({_jH_i},K).
\]
The adjunction yields for  $H_k$-modules
$M_k$, for $k\in\{i,j\}$, a natural isomorphism of vector spaces
\begin{equation} \label{eq:bimiso}
\Hom_{H_i}({_iH_j}\otimes_{H_j} M_j, M_i)\ra 
\Hom_{H_j}(M_j, {_jH_i}\otimes_{H_i} M_i), f\mapsto f^{\vee}.
\end{equation}

Quite similarly to
the representation theory of modulated graphs, in the sense of Dlab and
Ringel~\cite{DR1}, we have the following basic results 
from~\cite[Prp.~6.4]{qrel1} and~\cite[Prp.~7.1]{qrel1}.
\begin{Prop} \label{Prp:H-bas}
Set
$H:=H_K(C,D,\Ome)$. With 
$S:=\times_{i\in I} H_i$ we can consider
$\displaystyle B:=\bigoplus_{(i,j)\in\Ome} {_iH_j}$ as an $S$-$S$-bimodule and 
find:
\begin{itemize}
\item[(a)] 
$\displaystyle H\cong T_S(B):=\bigoplus_{j\in\NN} B^{\otimes_S j}$, 
i.e.~$H$ is a tensor algebra.
\item[(b)] There is a canonical short exact sequence of $H$-$H$-bimodules
\[
0\ra H\otimes_S B\otimes_S H \xrightarrow{\del} H\otimes_S H
\xrightarrow{\mathrm{mult}} H\ra 0,
\]
where $\del(h_l\otimes b\otimes h_r)=h_lb\otimes h_r - h_l\otimes bh_r$.
\end{itemize}
\end{Prop}
Note that the $H$-$H$-bimodules $H\otimes_SB\otimes_SH$ and $H\otimes_SH$ are
in general only projective as $H$-left- or right-modules, but not as bimodules. 
Anyway, the above sequence yields a functorial projective resolution
for certain modules which we are going to define now.
We say that a $H$-module $M$ is \emph{locally free} if $e_iM$ is a free
$H_i$-module for all $i\in I$. In this case we define
\[\rkv(M):= (\rk_{H_i} (e_iM))_{i\in I}.
\]
For example, there is a unique (indecomposable) locally free $H$-module
$E_i$ with $\rkv(E_i)=\alp_i$ for each $i\in I$. For later use we define
for all $\br\in\NN^I$ the module
$\bE^{\br}:=\oplus_{i\in I} E_i^{\br(i)}$, and observe that $\rkv(\bE^\br)=\br$.
Let us write down the following consequences of Proposition~\ref{Prp:H-bas}, 
see~\cite[Sec.~3.1]{qrel1} and~\cite[Cor.7.1]{qrel1}.
\begin{Cor} \label{Cor:H1}
For $H$ as above we have:
\begin{itemize}
\item[(a)] The projective and injective $H$-modules are locally free.
More precisely we have
\[
\rkv(He_{i_k})=\bet_{\bi,k}\quad\text{and}\quad\rkv(De_{i_k}H)=\gam_{\bi,k}\ 
\text{for } k\in I,
\]
where $\bi$ is a reduced expression for the Coxeter element of
$(C,\Ome)$.
\item[(b)]
Each locally free $H$-module $M$ has a functorial projective resolution
\[
0\ra H\otimes_SB\otimes_S M\xrightarrow{\del\otimes M} H\otimes_S M
\xrightarrow{\mathrm{mult}} M\ra 0.
\]
Moreover, if $M$ is not locally free, then $\pdim M=\infty$.
\item[(c)]
$H$ is 1-Iwanaga-Gorenstein, i.e. $\pdim({_HDH})\leq 1$ and 
$\idim({_HH})\leq 1$.
Moreover an $H$-module $M$ is locally free if and only if $\pdim(M)\leq 1$.
\end{itemize}
\end{Cor}

It follows that the Ringel (homological) bilinear form descends as the
non-symmetric bilinear form~\eqref{eq:non-symbil} to the Grothendieck
group of locally free modules, where we identify the classes of the
generalized simples $E_i$ with the  coordinate vector $\alp_i$ ($i\in I$),
see also~\cite[Prp.~4.1]{qrel1}.

\begin{Cor} \label{cor:ring-bil}
If $M$ and $N$ are locally free $H$-modules, we have
\[
\dim\Hom_H(M,N)-\dim\Ext^1_H(M,N)=\bil{\rkv(M),\rkv(N)}_{C,D,\Ome}.
\]
\end{Cor}
By combining~Corollary~\ref{Cor:H1} with 
standard results from Auslander-Reiten theory we obtain now the following
result.

\begin{Cor}\label{Cor:H2}
Let $M$ be an indecomposable, non projective, locally free $H$-module such 
that the  Auslander-Reiten translate $\tau_H M$ is locally free. Then
\[
\rkv(\tau_H M)=c\cdot(\rkv (M)),
\]
where $c=s_{i_n}\cdots s_{i_1}$ is the Coxeter element for $(C,\Ome)$.
Moreover, if we take $R\in\ZZ^{I\times I}$, such that $D\cdot R$ is the matrix
of $\bil{-,-}_{C,D,\Ome}$ with respect to the standard basis, 
we get $c=-R^{-1}(C-R)$.
\end{Cor}

This is the $K$-theoretic shadow of a  deeper connection between
the Auslander-Reiten translate and reflection functors, which we will
discuss in the next subsection.

\subsection{Auslander-Reiten theory and Coxeter functors}
By Proposition~\ref{Prp:H-bas} we may view $H=H_K(C,D,\Ome)$ as a tensor 
algebra. Thus, we identify a $H$-module $M$ naturally with a $S$-module 
$\bM=\oplus_{i\in I} M_i$ together with an  element $(M_{ij})_{(i,j)\in\Ome}$ of 
\begin{equation} \label{eqn:HM}
H(\bM):=\bigoplus_{(i,j)\in\Ome}
\in\Hom_{H_i}({_iH_j}\otimes_{H_j}M_j,H_i).
\end{equation}
Write $s_i(H):=H_k(C,D,s_i(\Ome))$ for any $i\in I$. If $k$ is
a sink of $Q^\circ(C,\Ome)$, we have for each $H$-module $M$ a canonical  
exact sequence
\begin{equation} \label{eq:M_in}
0\ra \Ker(M_{k,\inn})\ra \bigoplus_{j\in\Ome(k,-)} 
{_kH_j}\otimes_{H_j} M_j \xrightarrow{M_{k,\mathrm{in}}} M_k,
\text{ where } 
M_{k,\inn}=\oplus_{j\in\Ome(k,-)}M_{kj}.
\end{equation}
We can define now the BGP-reflection functor
\[
F_k^+\df\rep(H)\ra \rep(s_i(H)),\quad 
(F_k^+M)_i=\begin{cases} M_i &\text{if } i\neq k,\\
\Ker(M_{k,\inn}) &\text{if } i=k.\end{cases}
\]
We can moreover define in this situation dually the left adjoint
$F_k^-\df\rep(s_k(H))\ra\rep(H)$. Note that $k$ is a source of 
$Q^\circ(C,s_k\Ome)$. See~\cite[Sec.~9.2]{qrel1} for more details.
We observe  that the definitions imply easily the following:
\begin{Lem} \label{lem:refl1}
If $k$ is a sink for $\Ome$ and $M$ is a locally free
$H$-module which has no direct summand isomorphic to $E_k$ and $F_k^+(M)$
is locally free, then $\rkv(F_k^+M)=s_k(\rkv(M))$.
\end{Lem}

The proof of~\cite[Prp.~9.6]{qrel1} implies the following, less obvious result:
\begin{Lem} \label{lem:refl2}
Suppose that $k$ is a sink for $\Ome$ and $M$ a locally free rigid $H$-module, 
with no direct summand isomorphic to $E_k$, then $\Hom_H(M,E_k)=0$. 
\end{Lem}

We can interpret $F_k^+$ as a kind of APR-tilting functor~\cite{APR}. 
See~\cite[Sec.~9.3]{qrel1} for a proof of this non-trivial result.

\begin{Thm} \label{thm:apr}
Let $k$ be a sink of $Q^\circ(C,\Ome)$. 
Then $X:= {_HH}/He_k\oplus \tau^-He_k$ is a classical tilting module for
$H$. With $B:=\End_H(X)^\op$ we have an equivalence 
$S\df\rep(s_k(H))\ra\rep(B)$ such that the functors $S\circ F_l^+$ and
$\Hom_H(X,-)$ are isomorphic.
\end{Thm}

Standard tilting theory arguments and Auslander-Reiten theory, 
together with Lemma~\ref{lem:refl1} and Lemma~\ref{lem:refl2}
yield the following important consequence:
\begin{Cor} 
Let $k\in I$ be a sink for $\Ome$  and $M$ a locally free rigid
$H$-module, then $F_k^+(M)$ is a rigid, locally free $s_k(H)$-module.
\end{Cor}
Consider the algebra automorphism of $H$, which is defined by multiplying the
non-loop arrows of $Q(C,\Ome)$ by $-1$. It induces the so called
\emph{twist} automorphism $T\df\rep(H)\ra\rep(H)$.
Moreover, let $s_{i_n}\cdots s_{i_2}s_{i_1}$ be the Coxeter element
for $(C,\Ome)$, corresponding to the $+$-admissible sequence
$i_1, i_2,\ldots, i_n$, see Section~\ref{ssec:cox}. Now we can define the
\emph{Coxeter functor}
\[
C^+:= F_{i_n}^+\circ\cdots\circ F_{i_2}^+\circ F_{i_1}^+\df\rep(H)\ra\rep(H).
\]
Following ideas of P.~Gabriel and Ch.~Riedtmann~\cite[Sec.~5]{Gab}, 
by a careful comparison of
the definitions of the reflection functors and Auslander-Reiten translate,
we obtain the following result. See~\cite[Sec.~10]{qrel1} for the
lengthy proof.

\begin{Thm} \label{thm:GaRie}
With the $H$-$H$-bimodule $Y:=\Ext^1_H(DH,H)$ we have
an isomorphism of endofunctors of $\rep(H)$:
\[
\Hom_H(Y,-)\cong T\circ C^+
\]
If $M$ is locally free, we have  functorial isomorphisms 
\[
\tau_H(M)\cong \Hom_H(Y,M)\quad\text{and}\quad \tau_H^-M\cong Y\otimes_HM.
\]
In particular, in this case the Coxeter functor $C^+$ and the Auslander-Reiten
translate $\tau$ may be identified up to the twist $T$.
\end{Thm}

It is not true in general that the Auslander-Reiten translate of
a locally free $H$-module is again locally free. In~\cite[13.6-13.8]{qrel1}
several examples of this behavior are documented.
This motivates the following
definition. A  $H$-module $M$ is \emph{$\tau$-locally free}
if $\tau^k M$ is locally free for all $k\in\ZZ$. In particular,
rigid locally free modules are $\tau$-locally free. We call an indecomposable 
$H$-module \emph{preprojective}, resp.~\emph{preinjective}, if it is of
the form $\tau^{-k}(He_i)$ resp.~$\tau^k(De_iH)$ for some $k\in\NN_0$ and 
$i\in I$.
Thus, these modules are particular cases of rigid $\tau$-locally free modules.

\subsection{Dynkin type}
By combining the findings of previous section with standard Auslander-Reiten
theory and the characterization of Dynkin 
diagrams in Prop.~\ref{prp:Dynkin}, we obtain the following analog of 
Gabriel's theorem, see~\cite[Thm.~11.10]{qrel1}.

\begin{Thm} \label{Prp:H-fin}
Let $H=H_K(C,D,\Ome)$ be as above. There are only finitely
many isomorphism classes of indecomposable, $\tau$-locally free 
$H$-modules if and only if $C$ is of Dynkin type. In this case the map
$M\mapsto\rkv(M)$ induces a bijection between the isomorphism classes
of indecomposable, $\tau$-locally free modules and the positive roots
$\Del^+(C)$. Moreover, all these modules are preprojective and preinjective.
\end{Thm}

Note however, that even for $C$ of Dynkin type, the algebra $H(C,D,\Ome)$ 
is in most cases not of finite representation type, 
see~\cite[Prp.~13.1]{qrel1} for details. 

Let $C$ be a symmetrizable Cartan matrix of Dynkin type and 
$\bi=(i_1,i_2,\ldots,i_r)$ a reduced expression for the longest element
$w_0$ of the Weyl group $W$, which is $+$-admissible for the orientation $\Ome$.
With the notation of~\eqref{eq:rootlist} we abbreviate 
$\bet_j=\bet_{\bi,j}$ for
$j=1,\ldots, r$, and recall that this gives a complete list of the positive
roots. By Theorem~\ref{Prp:H-fin} we have for each $j$ a unique, locally free, 
indecomposable and rigid representation $M(\bet_j)$  
with $\rkv(M(\bet_j))=\bet_j$. 
\begin{Prop} \label{prp:homext}
With the above notations we have
\[
\bil{\bet_i,\bet_j}_{C,D,\Ome}=\begin{cases} 
\;\dim\Hom_H(M(\bet_i), M(\bet_j))  & \text{if } i\leq j,\\
-\dim\Ext^1_H(M(\bet_i), M(\bet_j)) &\text{if } i> j.\\
\end{cases}
\]
In particular,  $\Hom_H(M(\bet_i),M(\bet_j))=0$ if $i>j$ and 
$\Ext^1_H(M(\bet_i),M(\bet_j))=0$ if $i\leq j$.
\end{Prop}
In fact, by Theorem~\ref{thm:apr} and equation~\eqref{eq:refl-ring} we
may assume that either $i=1$ or $j=1$. In any case $M(\bet_1)=E_{i_1}$
is projective. In the first case we have $\Ext^1_H(E_1, M(\bet_j))=0$. 
In the second case we have $\Hom_H(M(\bet_i), E_{i_1})=0$ by 
Lemma~\ref{lem:refl2}. Now our claim follows by Corollary~\ref{cor:ring-bil}.

The next result is an easy adaptation of similar results by Dlab and 
Ringel~\cite{DR2} for species. The proof uses heavily 
Proposition~\ref{prp:homext}
and reflection functors. This version was worked out in Omlor's Masters
thesis~\cite{Oml}, see also~\cite[Sec.~5]{qrel5}.

\begin{Prop} \label{prp:nofilt}
With the same setup as above
let $k\in\{1,2,\ldots,r\}$ and $\bm=(m_1,\ldots, m_r)\in\NN^r$ such that 
$\bet_k=\sum_{j=1}^rm_j\bet_j$ and $m_k=0$. Then $M(\bet_k)$ admits
a non-trivial filtration by locally free submodules
\[
0=M_{(0)}\subset M_{(1)}\subset\cdots\subset M_{(r)}=M(\bet_k)
\]
such that $M_{(j)}/M_{(j-1)}\cong M(\bet_j)^{m_j}$ for $j=1,2,\ldots,r$.
It follows, that $M(\bet_k)$ has \underline{no} filtration by locally free 
submodules
\[
0=M^{(r)}\subset M^{(r-1)}\subset\cdots\subset M^{(0)}=M(\bet_k),
\]
such that $\rkv(M^{(j-1)}/M^{(j)})=m_j\bet_j$ for $j=1,2,\ldots, r$. 
\end{Prop}

\subsection{Generalized preprojective algebras}
Let $\bQ=\bQ(C)$ be the quiver which is obtained from $Q(C,\Ome)$ by inserting
for each $(i,j)\in\Ome$ additional $g_{ij}$ arrows 
$\alp_{ji}^{(1)},\ldots,\alp_{ji}^{(g_{ij})}$ from $i$ to $j$, and consider the
potential 
\[
W=\sum_{(i,j)\in\Ome}\sum_{k=1}^{g_{ij}}
(\alp_{ji}^{(k)}\alp_{ij}^{(k)}\eps_j^{-c_{ij}/g_{ij}}
-\alp_{ij}^{(k)}\alp_{ji}^{(k)}\eps_i^{-c_{ji}/g_{ij}}).
\]
The choice of $\Ome$ only affects the signs of the summands
of $W$. Recall that for a cyclic path $\alp_1\alp_2\cdots\alp_l$ in $\bQ$ 
by definition
\[
\partial_\alp^\cyc (\alp_1\alp_2\cdots\alp_l):=
\!\sum_{i\in\{j\in [1,l]\mid \alp_j=\alp\}}\!\!\!
\alp_{i+1}\alp_{i+2}\cdots\alp_l\alp_1\alp_2\cdots\alp_{i-1}.
\]
The generalized preprojective algebra of $H$ is
\[
\Pi=\Pi(Q,D):=K\bQ/
\bil{\partial^\cyc_\alp (W)\mid_{\alp\in\bQ_1},\ \eps_i^{c_i}\mid_{i\in I}}.
\]
It is easy to see that $\Pi$ does not depend on the choice
of $\Ome$, up to isomorphism. Notice that for $(i,j)\in\Ome$ we have
\[
\partial_{\ \alp_{ji}^{(k)}}^\cyc(W)=
\alp_{ij}^{(k)}\eps_j^{-c_{ij}/g_{ij}}-\eps_i^{-c_{ji}/g_{ij}}\alp_{ij}^{(k)}.
\]
It follows, that for any orientation $\Ome$ of $C$ we can equip $\Pi_K(C,D)$
with a $\NN_0$-grading by assigning each arrow $\alp_{ji}^{(k)}$ with 
$(i,j)\in\Ome$ degree $1$ and the remaining arrows get degree $0$. 
We write then
\[
\Pi_K(C,D)=\bigoplus_{i=0}^\infty  \Pi(C,D,\Ome)_i,
\]
and observe that $\Pi_K(C,D,\Ome)_0=H_K(C,D,\Ome)$. We obtain from
Theorem~\ref{thm:GaRie} the following
alternative description of our generalized preprojective algebra, which
justifies its name:

\begin{Prop} \label{Prp:Pi-name}
Let $C$ be a symmetrizable Cartan matrix with symmetrizer $D$, and $\Ome$
an orientation for $C$. Then, with $H=H_K(C,D,\Ome)$ we have
\[
\Pi(C,D,\Ome)_1\cong\Ext^1_H(DH,H)
\]
as an $H$-$H$-bimodule, moreover
\[
\Pi(C,D)\cong T_H(\Ext^1_H(DH,H))\quad\text{and}\quad
{_H\Pi(C,D)}\cong
\bigoplus_{i\in I, k\in\NN_0}
\tau_H^{-k}He_i.
\]
Here the first isomorphism is an isomorphism of $K$-algebras, and the
second one of $H$-modules.
\end{Prop}

Similarly to Proposition~\ref{Prp:H-bas} we have the following 
straightforward description of  our generalized preprojective
algebra as a tensor algebra modulo canonical 
relations~\cite[Prp.~6.1]{qrel1}, which yields a standard bimodule
resolution. See~\cite[Sec.~12.1]{qrel1} for the proof, where we 
closely follow~\cite[Lem.~3.1]{CBSh}. See also~\cite[Sec.~4]{BBK}.

\begin{Prop}\label{Prp:Pi-bimodres}
Let $C$ be a symmetrizable, connected Cartan matrix and $\Pi:=\Pi_K(C,D)$.
With $\oB:=\oplus_{(i,j)\in\Ome}({_iH_j}\oplus {_jH_i})$ we have 
$\Pi\cong T_S(\oB)/\bil{\partial_{\eps_i}^{\cyc}(W)\mid_{i\in I}}$,
where we interpret $\partial_{\eps_i}^\cyc(W)\in \oB\otimes_S\oB$ in
the obvious way. 
We obtain an exact sequence of $\Pi$-$\Pi$-bimodules
\begin{equation} \label{eq:PiRes}
\Pi\otimes_S\Pi\xrightarrow{f}\Pi\otimes_S\oB\otimes_S\Pi\xrightarrow{g}
\Pi\otimes_S\Pi\xrightarrow{h}\Pi\ra 0,
\end{equation}
where 
\[
f(e_i\otimes e_i)=
\partial_{\eps_i}^\cyc(W)\otimes e_i+e_i\otimes\partial^\cyc_{\eps_i}(W),\quad
g(e_i\otimes b\otimes e_j)=e_ib\otimes e_j-e_i\otimes be_j
\]
and $h$ is the multiplication map.
Moreover $\Ker(f)\cong\Hom_\Pi(D\Pi,\Pi)$ if $C$ is of Dynkin type, 
otherwise $f$ is injective.
\end{Prop}

We collect below several consequences, which can be found with detailed
proofs in~\cite[Sec.~12.2]{qrel1}. They illustrate that locally free
$\Pi$-modules behave in many aspects like modules over classical preprojective
algebras. Note that part~(b) is an extension
of Crawley-Boevey's remarkable formula~\cite[Lem.~1]{CB3}
\begin{Cor}\label{Cor:PiHomol}
Let $C$ be a connected, symmetrizable Cartan matrix, and $\Pi=\Pi_K(C,D)$
as above. Moreover, let $M$ and $N$ be locally free $\Pi$-modules.
\begin{itemize}
\item[(a)]
If $N$ finite-di\-mensional,
we have a functorial isomorphism
\[
\Ext^1_\Pi(M,N)\cong D\Ext^1_\Pi(N,M).
\]
\item[(b)]
If  $M$ and $N$ are finite-dimensional, we have
\[
\dim\Ext^1_\Pi(M,N)=\dim\Hom_\Pi(M,N)+\dim\Hom_\Pi(N,M)-(\rkv(M),\rkv(N))_{C,D}.
\]
\item[(c)]
If $C$ is not of Dynkin type, $\pdim(M)\leq 2$.
\item[(d)]
If $C$ is of Dynkin type, $\Pi$ is a finite-dimensional, self-injective
algebra and $\rep_\lfr(\Pi)$ is a 2-Calabi-Yau Frobenius category.
\end{itemize}
\end{Cor}
Similar to Cor.~\ref{Cor:H1}~(b) the  complex~\eqref{eq:PiRes} yields (the 
beginning of) a 
functorial projective resolution for all locally free $\Pi$-modules.
Thus (a), (b) and~(c) follow by exploring the symmetry of the above complex. 
For (d) we note that in this case $\Pi$ is finite-dimensional and ${_\Pi\Pi}$
is a locally free module by Thm.~\ref{Prp:H-fin} and Prp.~\ref{Prp:Pi-name}.

\section{Representation varieties}
\label{sec:rvar}
\subsection{Notation} \label{ssec:RV-not}
Let $K$ be now an algebraically closed field. For $Q$ a quiver and 
$\rho_j\in e_{t_j}(KQ_{\geq 2})e_{s_j}$ for $j=1,2,\ldots,l$ we set 
$A=KQ/\bil{\rho_1,\ldots,\rho_l}$. Note, that every finite dimensional basic
$K$-algebra is of this form. We abbreviate $Q_0=I$ and set 
for $\bd\in\NN_0^I$: 
\[
\Rep(KQ,\bd):=\times_{a\in Q_1}\Hom_K(K^{\bd(sa)},K^{\bd(ta)})\quad\text{and}\quad
\GL_\bd:=\times_{i\in I} \GL_{\bd(i)}(K).
\]
The reductive algebraic group $\GL_\bd$ acts on $\Rep(KQ,\bd)$ by conjugation, 
and the $\GL_\bd$-orbits correspond bijectively to the isoclasses of 
$K$-representations of $Q$. For $M\in\Rep(KQ,\bd)$ and 
$\rho\in e_iKQ e_j$ we can define  $M(\rho)\in\Hom_K(K^{\bd(j)},K^{\bd(i)})$ in
a natural way. We have then the $\GL_\bd$-stable, 
Zariski closed subset
\[
\Rep(A,\bd):=\{M\in\Rep(KQ,\bd)\mid M(\rho_i)=0\text{ for } j=1,2,\ldots, l\}.
\]
The $\GL_\bd$-orbits on $\Rep(A,\bd)$ correspond now to the isoclasses of 
representations of $A$ with dimension vector $\bd$.
It is in general a hopeless task to describe the irreducible components of
the affine variety $\Rep(A,\bd)$. 

\subsection{Varieties of locally free modules for $H$}
The set of locally free 
representations of $H=H_K(C,D,\Ome)$ is relatively easy to describe. Clearly,
for each locally free $M\in\rep(H)$ we have $\dimv(M)=D\cdot\rkv(M)$.
\begin{Prop} 
For $\br\in\NN^I$ we have the open subset
\[
\Rep_\lfr(H,\br):=\{M\in\rep(H,D\cdot\br)\mid M
\text{ is locally free}\}\subset \Rep(H,D\cdot\br),  
\]
which is irreducible and smooth with
$\dim\rep_\lfr(H,\br)=\dim\GL_{D\cdot\br} -\frac{1}{2}(\br,\br)_{C,D}$.
\end{Prop}
In fact, it is well known that the modules of projective dimension at most $1$
form always an open subset of $\rep(A,\bd)$. One verifies next
that $\Rep_\lfr(H,\br)$ is  a vector bundle over the $\GL_{D\cdot\br}$-orbit
$\cO(\oplus_{i\in I} E_i^{\br(i)})$, with the fibers isomorphic to the vector
space $H(\br):=H(\bE^\br)$, see~\eqref{eqn:HM}.

This yields the remaining claims. Note that the (usually) non-reductive
algebraic group
\[
G_\br:=\times_{i\in I} \GL_{\br(i)}(H_i)=\Aut_S(\oplus_{i\in I} E_i^{\br(i)})
\]
acts on the affine space $H(\br)$ naturally by conjugation, and 
the orbits are in bijection with isoclasses of locally free $H$-modules with
rank vector $\br$.

As a consequence, if $M$ and $N$ are rigid, locally free modules with
$\rkv(M)=\rkv(N)$, then already $M\cong N$, since the orbits of rigid modules
are open.

\subsection{Varieties of $\bE$-filtered modules for $\Pi$}
Recall the description of $\Pi_K(C,D)$ in Proposition~\ref{Prp:Pi-bimodres}. 
A $\oH:=T_S(\oB)$-module $M$ is given by a $S$-module 
$\bM=\oplus_{i\in I}M_i$  
such that $M_i$ is a $H_i$-module for $i\in I$, 
together with an element $(M_{ij})_{(i,j)\in\oOme}$ of 
\[
\oH(\bM):=\bigoplus_{(i,j)\in\oOme}\Hom_{H_i}({_iH_j}\otimes_{H_j} M_j, M_i),
\]
where $\oOme=\Ome\cap\Ome^{\op}$. Extending somewhat~\eqref{eq:M_in} we set
\begin{align*}
M_{i,\inn}:= &\left(\bigoplus_{j\in\oOme(i,-)} \sgn(i,j)M_{ij}\right)
\df \bigoplus_{j\in\oOme(i,-)}{_iH_j}\otimes_{H_j} M_j\ra M_i\quad\text{and}\\
M_{i,\out}:= &\left(\prod_{j\in\oOme(-,i)}M_{ji}^\vee\right)\df M_i\ra
\bigoplus_{j\in\oOme(-,i)} {_iM_j}\otimes_{H_j} M_j.
\end{align*}
We  define now for any $S$-module $\bM$, as above, the affine variety
\[
\Rep^{\fib}(\Pi,\bM):=\{(M_{ij})_{(i,j)\in\oOme}\in\oH(\bM)\mid
M_{k,\inn}\circ M_{k,\out}=0\text{ for all } k\in I\},
\]
and observe that the orbits of the, usually non-reductive, group  
$\Aut_S(\bM)$ on $\Rep^{\fib}(\Pi,\bM)$ correspond to the isoclasses of
possible structures of representations of $\Pi$ on $\bM$, since the
condition $M_{k,\inn}\circ M_{k,\out}$ corresponds to the relation
$\partial^\cyc_{\eps_k}(W)$.

Similarly to the previous section we can define the open subset
\[
\Rep_\lfr(\Pi,\br):=\{M\in\Rep(\Pi,D\cdot\br)\mid M\text{ locally free}\}
\subset\Rep(\Pi, D\cdot\br),
\]
and observe that $\Rep_\lfr(\Pi,\br)$ is a fiber bundle over the
$\GL_{D\cdot\br}$-orbit $\cO(E^{\br})$, with typical fiber 
$\Rep^{\fib}(\Pi,\bE^{\br})$. Finally we define for any projective $S$-module
$\bM$ the constructible subset
\[
\Pi(\bM)=\{ (M_{ij})_{(i,j)\in\oOme}\in\Rep^\fib(\Pi,\bM)\mid ((M_{ij})_{ij},\bM)
\text{ is } \bE \text{-filtered}\}.
\]
Here, a $\Pi$-module $X$ is $\bE$-filtered if it admits a flag of submodules
$0=X_{(0)}\subset X_{(1)}\subset\cdots\subset X_{(l)}=X$, such that for all $k$
we have $X_{(k)}/X_{(k-1)}\cong E_{i_k}$ for some $i_1,i_2,\ldots,i_l\in I$.
Note that for $C$ symmetric and $D$ trivial this specializes to Lusztig's
notion of a nilpotent representation for the preprojective algebra of
a quiver. However, if $C$ is not symmetric even in the Dynkin case 
there exist finite-dimensional, 
locally free $\Pi$-modules which are not $\bE$-filtered, 
see~\cite[Sec.~8.2.2]{qrel4} for an example.

We consider $\Pi(\br)$ with the Zariski topology and call it by a slight
abuse of notation a variety. In any case, it makes sense to speak of the 
dimension of $\Pi(\br)$ and we can consider the set
\[
\Irr(\Pi(\br))^{\max}
\]
of top-dimensional irreducible components of $\Pi(\br)$. 
 
\begin{Thm} \label{Prp:IrrPi}
Let $C$ be a symmetrizable generalized Cartan matrix with symmetrizer $D$
and $H=H_K(C,D,\Ome), \Pi=\Pi_K(C,D)$ for an algebraically closed field $K$.
For the spaces $\Pi(\br)$ of $\bE$-filtered representations of $\Pi$ we have 
\begin{itemize}  
\item[(a)]
$\displaystyle
\dim \Pi(\br)=\dim H(\br)=\sum_{i\in I}c_i\br(i)^2-\frac{1}{2}(\br,\br)_{C,D}
$
for all $\br\in\NN^I$.
\item[(b)] The set $\cB:=\coprod_{\br\in\NN^I} \Irr(\Pi(\br))^{\max}$ has a natural
structure of a crystal of type $B_C(-\infty)$ in the sense of Kashiwara.
In particular, we have 
\[
\abs{\Irr(\Pi(\br))^{\max}}=\dim U(\fn)_\br,
\]
where $U(\fn)$  the universal
enveloping algebra of the positive part $\fn$ of the Kac-Moody Lie algebra
$\fg(C)$.  
\end{itemize}
\end{Thm}
We will sketch in the next two sections a proof of these two statements, 
which are the main result of~\cite{qrel4}.

\subsection{Bundle constructions}
The  bundle construction in this section is crucial. 
It is our version~\cite[Sec.~3]{qrel4} of 
Lusztig's construction~\cite[Sec.~12]{L1}.

For $m\in\NN$ we denote by $\cP_m$ the set of  sequences
of integers $\bp=(p_1, p_2,\ldots,p_t)$ with 
$m\geq p_1\geq p_2\geq\cdots\geq p_t\geq 0$. Obviously $\cP_{c_k}$ parametrizes
the isoclasses of $H_k$-modules, and we define 
$H_k^\bp=\oplus_{j=1}^t H_k/(\eps_k^{p_j})$. 
For $k\in I$ and $M\in\rep(\Pi)$ we set 
\[
\fac_k(M):= M_k/\Ima(M_{k,\inn})\quad\text{and}\quad \sub_k(M):=\Ker(M_{k,\out}).
\]
With this we can  define  
\begin{align*}
\Pi(\bM)^{k,\bp}&=\{M\in\Pi(\bM)\mid \fac_k(M)\cong H_k^{\bp}\} \text{ and }\\
\Pi(\bM)_{k,\bp}&=\{M\in\Pi(\bM)\mid \sub_k(M)\cong H_k^{\bp}\} 
\end{align*}
for $\bp\in\cP_{c_k}$.
We abbreviate $\Pi(\bM)^{k,m}=\Pi(\bM)^{k,c_k^m}$. In what follows, we will focus
our exposition on the varieties of the form $\Pi(\bM)^{k,\bp}$, however one
should be aware that similar statements and constructions hold for the
dual versions $\Pi(\bM)_{k,\bp}$.

For an $\bE$-filtered
representation $M\in\rep(\Pi)$ there exists always a $k\in I$ such that
$\fac_k(M)$, viewed as an $H_k$-module, has a non-trivial free summand. 
It is also important to observe that $\Pi(\bM)^{k,0}$ is an open subset of
$\Pi(\bM)$.

Fix now $k\in I$, let $\bM$ be a projective $S$-module and $\bU$ be a 
proper, projective  $S$-submodule of $\bM$ with $U_j=M_j$ for all $j\neq k$.
Thus, $\bM/\bU\cong E_k^r$ for some $r\in\NN_+$, and can choose a (free) 
complement $T_k$, such that $M_k=U_k\oplus T_k$.
For partitions $\bp=(c_k^r,q_1, q_2,\ldots,q_t)$ and $\bq=(q_1,\ldots, q_t)$
in $\cP_{c_k}$ we set moreover 
$\Hom^{\mathrm{inj}}_S(\bU,\bM):=\{f\in\Hom_S(U,M)\mid f\text{ injective}\}$,
and define
\[
Y^{k,\bp,\bq}:=\{(U,M,f)\in\Pi(\bU)^{k,\bq}\!\!\times\Rep^\fib(\Pi,\bM)\times
\Hom_S^{\mathrm{inj}}(\bU,\bM)\mid
f\in\Hom_\Pi(U,M)\}.
\]
Note that for $(U,M,f)\in Y^{k,\bp,\bq}$ we have in fact $M\in\Pi(\bM)^{k,\bp}$, 
and that the group $\Aut_S(\bU)$ acts freely on $Y^{k,\bp,\bq}$ via 
\[
g\cdot (U,M,f):= ((g_iU_{ij}(\id\otimes g_j^{-1}))_{(i,j)\in\oOme},\ M, 
\ g\cdot f^{-1}).
\]
\begin{Lem}\label{Lem:genLusztig}
Consider in the above situation the diagram
\[
\xymatrix{
&\ar[ld]_{p'} Y^{k,\bp,\bq}\ar[rd]^{p''}\\
\Pi(\bU)^{k,\bq}\times\Hom^{\mathrm{inj}}_S(\bU,\bM) & &\Pi(\bM)^{k,\bp}}
\]
with $p'(U,M,f)=(U,f)$ and $p''(U,M,f)=M$. Then the following holds:
\begin{itemize}
\item[(a)] $p'$ is a vector bundle of rank $m$, where
\[
m=\sum_{j\in\oOme(-,k)}\dim_K \Hom_K(T_k,{_kM_j}\otimes_{H_j}M_j) -
\dim_K\Hom_{H_k}(T_k,\Ima(U_{k,\inn})).
\]
\item[(b)]
$p''$ is a fiber bundle with smooth irreducible fibers isomorphic to
\[
\Aut_S(\bU)\times \Gr^{T_k}_{H_k}(H_k^{\bp}),
\]
where $\Gr^{T_k}_{H_k}(H_k^{\bp}):=\Hom_{H_k}^{\mathrm{surj}}(H_k^{\bp},T_k)/\Aut_{H_k}(T_k)$.
\end{itemize}
\end{Lem}

\begin{Cor}
In the situation of Lemma~\ref{Lem:genLusztig}, the correspondence
\[
Z' \mapsto p''({p'}^{-1}(Z'\times \Hom_S^{\mathrm{inj}}(\bU,\bM)):=Z''
\]
induces a bijection between the sets of irreducible components 
$\Irr(\Pi(U)^{k,\bq})$ and $\Irr(\Pi(\bM)^{k,\bp})$.
Moreover we have then
\[
\dim Z''-\dim Z'=\dim H(\bM)-\dim(\bU).
\]
\end{Cor} 
Note, that this implies already part~(a) of Theorem~\ref{Prp:IrrPi}. 
In fact the Corollary allows us to conclude by induction that
$\dim\Pi(\br)\leq\dim\Rep^{\fib}(H,\br)$. On the other hand, we can
identify $H(\br)$ with an irreducible component of $\Pi(\br)$.

\subsection{Crystals}
For $M\in\rep(\Pi)$ and $j\in I$ there are two canonical short exact sequences
\[
0\ra K_j(M)\ra  M\ra \fac_j(M)\ra 0\quad\text{ and }\quad
0\ra \sub_j(M)\ra M \ra C_j(M)\ra 0.
\]
We define recursively that $M$ is a \emph{crystal module} if 
$\fac_j(M)$ and $\sub_j(M)$ are locally free for all $j\in I$, and
$K_j(M)$ as well as $C_j(M)$ are crystal modules for all $j\in I$.
Clearly, if $M$ is a crystal module, for all $j\in I$ there exist
$\vph_j(M),\vph_j^*(M)\in\NN$ such that 
\begin{equation} \label{eq:phi}
\sub_j(M)\cong E_j^{\vph_j(M)} \text{and }\
\fac_j(M)\cong E_j^{\vph_j^*(M)}. 
\end{equation}
Note moreover, that crystal modules
are by construction $\bE$-filtered. It is now easy to see that for
all projective $S$-modules $\bM$ the set
\[
\Pi(\bM)^\crys:=\{M\in\Pi(\bM)\mid M\text{ is a crystal representation}\}
\]
is a constructible subset of $\Pi(\bM)$. The following result 
from~\cite[Sec.~4]{qrel4} is crucial for
the proof of Proposition~\ref{Prp:IrrPi}~(b). It has no counterpart for
the case of trivial symmetrizers.
\begin{Prop} \label{prp:cryst}
For each projective $S$-module $\bM$ the set $\Pi(\bM)^\crys$ is a dense and
equidimensional subset of the union of all top dimensional irreducible 
components of $\Pi(\bM)$.
\end{Prop}

This allows us in particular to define for all $Z\in\Irr(\Pi(\bM))^{\max}$ 
and $i\in I$ the value $\vph_i(Z)$, see~\eqref{eq:phi}, 
such that for a dense open subset
$U\subset Z$ we have $\vph_i(M)=\vph_i(Z)$ for all $M\in U$. Similarly,
we can define $\vph_i^*(Z)$.

Next we set
\[
\Irr(\Pi(\br)^{i,p})^{\max}:=
\{Z\in\Irr(\Pi(\br)^{i,p})\mid\dim Z=\dim H(\br)\}
\]
for $i\in I$ and $p\in\NN_0$, and similarly $\Irr(\Pi(\br)_{i,p}$. By
Lemma~\ref{Lem:genLusztig} we get a bijection
\[
e_i^*(\br,p)\df\Irr(\Pi(\br)^{i,p})^{\max}\ra
\Irr(\Pi(\br+\alp_i)^{i,p+1})^{\max}, Z\mapsto p''(p'^{-1}(Z\times J_0))
\]
Similarly we obtain a bijection
\[
e_i(\br,p)\df\Irr(\Pi(\br)_{i,p})^{\max}\ra\Irr(\Pi(\br+\alp_i)_{i,p+1})^{\max}.
\]
This allows us to define for all $\br\in\NN^I$ the operators
\[
\tilde{e}_i\df\Irr(\Pi(\br))^{\max}\ra\Irr(\Pi(\br+\alp_i)),
Z\mapsto \overline{e_i(\br,\vph_i(Z))(Z^\circ)},
\]
where $Z^\circ\in\Irr(\Pi(\br)_{i,\vph_i(Z)})^{\max}$ is the unique irreducible 
component with $\overline{Z^\circ}=Z$. 
Similarly, we can define the operators $\tilde{e}^*_i$ 
in terms of the bijections $e_i^*(\br,p)$. 
We define now 
\begin{equation} \label{eq:cB}
\cB:= \coprod_{\br\in\NN_0} \Irr(\Pi(\br))^{\max} 
\text{ and } \wt\df\cB\ra\ZZ^I,\ Z\mapsto\rkv(Z).
\end{equation}
It is easy to see that 
$(\cB,\wt, (\tilde{e}_i,\vph_i)_{i\in I})$ is special case of a lowest
weight crystal in the sense of Kashiwara~\cite[Sec.~7.2]{K1}, namely we have
\begin{itemize}
\item 
$\vph_i(\tilde{e}_i(b))=\vph_i(b)+1,\ \wt(\tilde{e}_i(b))=\wt(b)+\alp_i$,
\item 
with $\{b_-\}:=\Irr(\Pi(0))^{\max}$, for each $b\in\cB$ there exists
a sequence $i_1,\ldots,i_l$ of elements of $I$ with 
$\tilde{e}_{i_1}\tilde{e}_{i_2}\cdots\tilde{e}_{i_l}(b_-)=b$,
\item
$\vph_i(b)=0$ implies $b\not\in\Ima(\tilde{e}_i)$.
\end{itemize}
Similarly $(\cB,\wt, (\tilde{e}_i^*,\vph_i^*)_{i\in I})$ is
a lowest weight crystal with the same lowest weight element $b_-$.

\begin{Lem} The above defined operators and functions on $\cB$ fulfill
additionally the following conditions:
\begin{itemize}
\item[(a)]
If $i \not= j$, then
$\tilde{e}_i^*\tilde{e}_j(b) = \tilde{e}_j\tilde{e}_i^*(b)$.
\item[(b)]
For all $b \in B$ we have $\vph_i(b) + \vph_i^*(b) - 
\bil{\wt(b),\alpha_i} \ge 0$.
\item[(c)]
If $\vph_i(b) + \vph_i^*(b) -  \bil{\wt(b),\alp_i} = 0$, 
then $\tilde{e}_i(b) = \tilde{e}_i^*(b)$.
\item[(d)]
If $\vph_i(b) + \vph_i^*(b) - \bil{\wt(b),\alp_i} \ge 1$, 
then $\vph_i(\tilde{e}_i^*(b)) = \vph_i(b)$ and \\
${\vph_i^*(\tilde{e}_i(b)) = \vph_i^*(b)}$.
\item[(e)]
If $\vph_i(b) + \vph_i^*(b) - \bil{\wt(b),\alp_i} \ge 2$, 
then $\tilde{e}_i\tilde{e}_i^*(b) = \tilde{e}_i^*\tilde{e}_i(b)$.
\end{itemize}
\end{Lem}
The proof of this Lemma in~\cite[Sec.~5.6]{qrel4} uses the homological 
features of locally free $\Pi$-modules from Corollary~\ref{Cor:PiHomol}
in an essential way. 
Note that here, by definition, $\bil{\br,\alp_i}=(C\cdot\br)_i$.

Altogether this means, by a criterion of Kashiwara and 
Saito~\cite[Prp.~3.2.3]{KS}, which we use here in a reformulation due
to Tingley and Webster~\cite[Prp.~1.4]{TW}, that
$(\cB,\wt,(\tilde{e}_i,\vph_i)_{i\in I}) \cong 
(\cB,\wt,(\tilde{e}_i^*,\vph_i^*)_{i\in I}) \cong B_C(-\infty)$.
Here, $B_C(-\infty)$ is the  crystal graph of the quantum group $U_q(\fn(C))$. 
This implies part~(b) of Theorem~\ref{Prp:IrrPi}.

\begin{Rem} We did not give here Kashiwara's general definition of a 
crystal graph, or that of a lowest weight crystal associated to a dominant 
integral weight.
The reason is that, due to limitations of space, we  can not to set up 
the, somehow unwieldy, notations for the integral weights of a Kac-Moody
Lie algebra. The interested reader can look up the relevant definitions,
in a form which is compatible with these notes, in~\cite[Sec.~5.1, 5.2]{qrel4}.
\end{Rem}

\section{Algebras of constructible functions}
\label{sec:const}
\subsection{Constructible functions and Euler characteristic}
Recall that  the topological Euler characteristic, defined in terms of singular 
cohomology with compact support and rational coefficients, defines a ring 
homomorphism from the Grothendieck ring of complex varieties to the integers. 

By definition, a constructible function $f\df X\ra\CC$ on a complex algebraic
variety $X$ has finite image, and $f^{-1}(c)\subset C$ is a constructible set
for all $c\in\CC$. By the above remark it makes sense to define
\[
\int_{x\in X} f d\chi := \sum_{c\in\CC} c\,\chi(f^{-1}(c)).
\]
If $\vph\df X\ra Y$ is a morphism of varieties, we can define the 
\emph{push forward} of constructible functions via 
$(\vph_*(f))(y):=\int_{x\in\vph^{-1}(y)}f d\chi$. This is functorial in the sense 
that $(\psi\circ\vph)_*(f)=\psi_*(\vph_*(f))$ for $\psi\df Y\ra Z$ an other
morphism, by result of McPherson~\cite[Prp.~1]{MPh}. See also~\cite[Sec.~3]{Jo}
for a careful discussion.

\subsection{Convolution algebras as enveloping algebras}
Let $A=\CC Q$ as in Section~\ref{ssec:RV-not}.  We consider for a dimension 
vector
$\bd\in \NN^I$ the vector space $\cF(A)_\bd$ of constructible functions 
$f\df\Rep(A,\bd)\ra\CC$ which are constant on $\GL_\bd(\CC)$-orbits and
set
\[
\cF(A):= \bigoplus_{\bd\in\NN^i} \cF(A)_\bd.
\]
Following Lusztig~\cite{L1}
$\cF(A)$ has the structure of a unitary, graded associative algebra.
The multiplication  is defined by
\[
(f*g)(X)= \int_{U\in\Gr^A_\bd(X)}f(U)g(X/U) d\!\chi,
\]
where $f\in\cF(A)_\bd$, $g\in\cF(A)_\be$, $X\in\Rep(A,\bd+\be)$, and
$\Gr_\bd^A(X)$ denotes the quiver Grassmannian of $\bd$-dimensional 
subrepresentations of $X$. The associativity of the multiplication
follows easily from the functoriality of the push-forward of constructible
functions. We have an algebra homomorphism
\begin{equation} \label{eq:comult}
c\df \cF(A) \ra \cF(A\times A),\ \text{with } (c(f))(X,Y)=f(X\oplus Y), 
\end{equation}
see for example~\cite[Sec.~4.3]{qrel3}. The proof depends crucially on the
Białynicki-Birula result about the fixpoints of algebraic torus 
actions~\cite[Cor.~2]{BB}. 
This fails for example over the real numbers.

\begin{Rem} \label{rem:M-filt}
If $\bX=(X_j)_{j\in J}$ is a family of indecomposable representations of $A$,
we  define the characteristic functions $\vth_j\in\cF_{\dimv X_j}(A)$
of the $\GL_{\dimv X_j}$-orbit $\cO(X_j)\subset\Rep(A,\dimv X_j)$ and 
consider the graded subalgebra $\cM(A)=\cM_\bX(A)$ of $\cF(A)$, which is
generated by the $\vth_j$. Clearly, the homogeneous components of $\cM$
are finite dimensional. If $\bj=(j_1,j_2,\ldots,j_l)$ is a sequence of elements
of $j$ we have by the definition of the multiplication
\[
\vth_{j_1}*\vth_{j_l}*\cdots *\vth_{j_l}(X)=\chi(\Fl_{\bX,\bj}^A(M)),
\]
where $\Fl_{\bX,\bj}^A(M)$ denotes the variety of all flags of submodules
\[
0=M^{(0)}\subset M^{(1)}\subset\cdots\subset M^{(l)}=M
\]
with $M^{(k)}/M^{(k-1)}\cong X_{j_k}$ for $k=1,2,\ldots,l$. In particular,
if $M$ has no filtration with all factors isomorphic to some $X_j$, we have
$f(M)=0$ for all $f\in\cM(A)_{\dimv M}$. See~\cite[Lemma~4.2]{qrel3}.
\end{Rem}

\begin{Lem} The morphism $c$ from~\eqref{eq:comult} induces a 
comultiplication
$\Del\df \cM(A)\ra \cM(A)\otimes\cM(A)$ with 
$\Del(\vth_j)=\vth_j\otimes 1+ 1\otimes\vth_j$ for all $j$. With this structure
$\cM$ is a cocommutative Hopf algebra, which is isomorphic to the enveloping 
algebra $U(\cP(\cM))$ of the Lie algebra of its primitive elements
$\cP(\cM)$.
\end{Lem}
See~\cite[Prp.~4.5]{qrel3} for a proof. Recall, that an element $x$ of a 
Hopf algebra is called \emph{primitive} iff $\Del(x)=x\otimes 1+1\otimes x$.
It is straightforward to check that the primitive elements of a Hopf algebra
form a Lie algebra under the usual commutator $[x,y]=xy-yx$. 

\begin{Rem} \label{rem:primsupp}
It is important to observe that, by the definition of the comultiplication,
the support of any primitive element of $\cM$ consists of indecomposable,
$\bX$-filtered modules. In fact, for $f\in\cP(\cM)$ and $M,N\in\rep(A)$
we have $f(M\oplus N)=cf(M,N)=(f\otimes 1+ 1\otimes f)(M,N)$. 
See~\cite[Lem.~4.6]{qrel3}.
\end{Rem}

We are here interested in the two special cases when $A=H_\CC(C,D,\Ome)$ or
$A=\Pi_\CC(C,D)$ and $\bX=\bE=(E_i)_{i\in I}$. Note that by 
Remark~\ref{rem:M-filt} only locally free modules can appear
in the support of any $f\in\cM_\bE(H)$. Similarly, the support of any
$f\in\cM_\bE(\Pi)$ consists only of $\bE$-filtered modules.
For this reason we will consider in what follows, both $\cM_\bE(H)$ and
$\cM_\bE(\Pi)$ as graded by \emph{rank vectors}. In other words, from now on
\[
\cM(H):=\cM_\bE(H)=\bigoplus_{\br\in\NN^I} \cM_\br(H)\quad\text{and}\quad\
\cM(\Pi)=\cM_\bE(\Pi)=\bigoplus_{\br\in\NN^I} \cM_\br(\Pi),
\]
where we may consider the the elements of $\cM_\br(H):=\cM_\bE(H)_{D\cdot\br}$ as 
constructible functions on $H(\br)$. 
Similarly we may consider the elements of $\cM_\br(\Pi):=\cM_\bE(\Pi)_{D\cdot\br}$
as constructible functions on $\Pi(\br)$. 

\subsection{$\cM_\bE(H)$ and a dual PBW-basis   in the Dynkin case} \label{ssec:pbw}
We have the following basic result from~\cite[Cor.~4.10]{qrel3}. 
\begin{Prop} \label{prp:pbw1}
Let $C$ be a symmetrizable Cartan matrix, $D$ a symmetrizer and $\Ome$ an
orientation for $C$. With $H=H_\CC(C,D,\Ome)$ we have an surjective Hopf
algebra homomorphism
\[
\eta_H\df U(\fn(C))\ra \cM_\bE(H)\ \text{defined by}\quad 
e_i\mapsto \vth_i (i\in I).
\]
\end{Prop}
The main point is to show that for the $\vth_i$ $(i\in I)$ fulfill the 
Serre relations~\eqref{eq:SerRel}. More precisely we need that the
primitive elements
\[
\vth_{ij}:=(\ad \vth_i)^{1-c_{ij}}(\vth_j)\in\cP(\cM(H))_{(1-c_{ij})\alp_i+\alp_j}
\quad (i\neq j)
\]
actually vanish. For this it is enough, by Remark~\ref{rem:primsupp}, 
to show that there exists
no \emph{indecomposable}, locally free $H$-module $M$ with 
$\rkv(M)=(1-c_{ij})\alp_i+\alp_j$. This is carried out in the proof 
of~\cite[Prp.~4.9]{qrel3}.

The proof of the following result, which is~\cite[Thm.~6.1]{qrel3}, occupies
the major part of that paper.
\begin{Prop} 
Let $C$ be a symmetrizable Cartan matrix of Dynkin type,
$D$ a symmetrizer and $\Ome$ an orientation for $C$ and $H=H_\CC(C,D,\Ome)$.
Then for each  positive root $\bet\in\Del^+$ there exists a primitive element
$\vth_\bet\in\cP(\cM(H))_\bet$ with $\vth_\bet(M(\bet))=1$.
\end{Prop} 

The idea of the proof is as follows: By~\cite[Cor.~1.3]{qrel2} for any
$\bet\in\Del^+(C)$ and any sequence $\bi$ in $I$, the
Euler characteristic $\chi(\Fl^H_{\bE,\bi}(M(\bet)))$ is independent of
the choice of the symmetrizer $D$. So, we may assume that $C$ is connected
and $D$ minimal. In the symmetric (quiver) case, our claim follows now
by Schofield's result~\cite{S}, who showed that in this case 
$\cP(\cM(H))$ can be identified with $\fn(C)$. By Gabriel's theorem in
this case the $\vth_\bet$ are the characteristic function of the 
$\GL_\bet$-orbit of $M(\bet)$.

In the remaining cases,
we construct the $\vth_\bet$ by induction on the height of $\beta$ in
terms of (iterated) commutators of ``smaller'' $\vth_\gam$. Note however
that in this case this construction is delicate since the support of the
$\vth_\bet$ may contain several indecomposable, locally free modules.
See for example~\cite[Sec.~13.2(d)]{qrel5}.
\medskip

Since in the Dynkin case all weight spaces of $\fn(C)$ are one-dimensional,
the main result of~\cite{qrel3}, Theorem~1.1~(ii), follows easily:
\begin{Thm}
If $C$ is of Dynkin type, the Hopf algebra homomorphism $\eta_H$ is
an isomorphism.
\end{Thm}

Recall the notation used in Proposition~\ref{prp:homext}. In particular,
$\bi$ is a reduced expression for the longest element $w_0\in W(C)$, which
is $+$-adapted to $\Ome$, and $\bet_k=\bet_{\bi,k}$ for $k=1,2,\ldots, r$. 
Let us abbreviate
\[
\vth_\bm := \frac{1}{m_r!\cdots m_1!}
\vth_{\bet_r}^{m_r}*\cdots *\vth_{\bet_1}^{m_1}\quad\text{and}\quad
M(\bm) := \oplus_{k=1}^r M(\bet_k)^{m_k}
\]
for $\bm=(m_1, m_2,\ldots, m_r)\in\NN^r$.
By the above results $(\vth_\bm)_{\bm\in\NN^r}$ is a normalized PBW-basis of 
$\cM(H)\cong U(\fn(C))$ in the Dynkin case.

Moreover we consider the graded dual $\cM(H)^*$ of $\cM(H)$, and the evaluation
form $\del_{M(\bm)}\in\cM(H)^*$ with $\del_\bm(f):=f(M(\bm))$.  By the 
definition of the comultiplication in $\cM(H)$, the graded dual is a 
commutative Hopf algebra, and $\del_{M(\bm)}\cdot\del_{M(\bn)}=\del_{M(\bm+\bn)}$. 
Our next result is essentially~\cite[Thm.~1.3]{qrel5}.

\begin{Prop} \label{prp:dual-PBW}
With the above notation we have
\[
\del_{M(\bm)}(\vth_{\bn})=\del_{\bm,\bn} \text{ for all } \bm,\bn\in\NN^r.
\]
Thus $(\del_{M(\bm)})_{\bm\in\NN^r}$ is a basis of $\cM(H)^*$ which is dual to
the PBW-basis $(\vth_\bm)_{\bm\in\NN^r}$, and
$\cM(H)^*=\CC[\del_{M(\bet_1)},\ldots,\del_{M(\bet_r)}]$.
\end{Prop}
In the quiver case (with trivial symmetrizer) this result is easy to prove,
since with Gabriel's theorem and Proposition~\ref{prp:homext} follows quickly
that $\vth_{\br}$ is the characteristic function of the orbit of $M(\br)$.
However, in our more general setting, already the $\vth_{\bet_k}$ are usually
not the characteristic function of $M(\bet_k)$, as we observed above.
The more sophisticated Proposition~\ref{prp:nofilt} implies, 
by the definition of the multiplication in $\cM_\bE(H)$, that
$\vth_{\bm}(M(\bet_k))=0$  if $\bm\neq\be_k$, the $k$-th unit vector.
The remaining claims follow now by formal arguments, 
see the proof of~\cite[Thm.~6.1]{qrel5}.

\begin{Rem} For $M\in\rep_\lfr(H)$ and $\be\in\NN^I$ we 
we introduce the quasi-projective variety
\[
\Grlf_\be^H(M):=\{U\subset M\mid U\text{ locally free submodule and } \rkv(U)=\be\},
\]
which is an open subset of the usual quiver Grassmannian $\Gr_{D\cdot\be}^H(M)$.
With this notation we can define
\[
F_M:=\sum_{\be\in\NN^I}\chi(\Grlf_\be^H(M))Y^\be\in\ZZ[Y_1,\ldots,Y_n]
\text{ and }
\bg_M:= -R\cdot\rkv(M),
\]
where $R$ is the matrix introduced in Corollary~\ref{Cor:H2}. By the main result
of~\cite{qrel5} this yields for $M=M(\bet)$ with $\bet\in\Del^+(C)$ 
the $F$-polynomial and $g$-vector, in the sense of~\cite{FZ3}, for all cluster
variables of a finite type cluster algebra~\cite{FZ2} of type $C$ with respect
to an acyclic seed defined by $\Ome$. The proof is based on
Proposition~\ref{prp:dual-PBW}, and on the description by Yang and 
Zelevinsky~\cite{YZ} of the $F$-polynomial of a cluster variable in terms
of generalized minors. 
\end{Rem}

\subsection{Semicanonical functions and the support 
conjecture for $\cM_\bE(\Pi)$}
Recall, that we abbreviate $\Pi=\Pi_\CC(C,D)$ for a symmetrizable Cartan matrix
$C$ with symmetrizer $D$.
By definition $\cM(\Pi)=\cM_\bE(\Pi)\subset\cF(\Pi)$ is generated by the 
functions $\tth_i\in\cM_{\alp_i}(\Pi)$ for $i\in I$, 
where $\tth_i$ is the characteristic
function of the orbit of $E_i$, viewed as a $\Pi$-module.
We use here the notation $\tth_i$ rather
than $\vth_i$ to remind us that the multiplication is now defined in 
terms of constructible functions on a larger space.
More precisely, we have for each $\br\in\NN^I$ an injective 
$\Aut_S(\bE^\br)$-equivariant, injective morphism of varieties
\[
\iota_\br\df H(\br)\ra\Pi(\bd).
\]
These morphisms induce, via restriction, a surjective morphism of 
graded Hopf algebras
\[
\iota_\Ome^*\df\cM(\Pi)\ra\cM(H),\quad \tth_i\mapsto\vth_i\ \text{for } i\in I. 
\]
The proof of the following result is, almost verbatim,
the same induction argument as the one used by Lusztig~\cite{L2}, 
see~\cite[Lem.~7.1]{qrel4}.

\begin{Lem} \label{lem:semican}
Let $\br\in\NN^I$. For each $Z\in\Irr(\Pi(\br))^{\max}$ there
exists an open dense subset $U_Z\subset Z$ and a function $f_Z\in\cM_\br(\Pi)$
such that for $Z,Z'\in\Irr(\Pi(\br))^{\max}$ and any $u'\in U_{Z'}$ we have
\[
f_Z(u')=\del_{Z,Z'}.
\]
In particular, the functions $(f_Z)_{Z\in\Irr(\Pi(\br))^{\max}}$ are linearly 
independent in $\cM_\br(\Pi)$.
\end{Lem}
Note however, that the result is not trivial since we claim that
the $f_Z\in\cM_\be(\Pi)$ and not in the much bigger space 
$\cF(\Pi)_{C\cdot\br}$.
On the other hand, it is important to observe
that the inductive construction of the \emph{semicanonical functions} 
$f_Z$ involves some choices.

As in Section~\ref{ssec:pbw}, we define now for each $i\neq j$ in $I$ 
the primitive element
\[
\tth_{ij}=(\ad\tth_i)^{1-c_{ij}}(\tth_j)\in\cP(\cM(\Pi)).
\]
Unfortunately, we have the following result, which is a combination of 
Lemma~6.1, Proposition~6.2 and Lemma~6.3 from~\cite{qrel4}.
\begin{Lem} \label{lem:noSerre}
Suppose with the above notations that  $c_{ij}<0$.
\begin{itemize}
\item[(a)] 
If $c_i\geq 2$ then there exists an indecomposable,
$\Pi=\Pi(C,D)$-module $X=X_{(ij)}$ with
$\rkv(X_{(ij)}=(1-c_{ij})\alp_i+\alp_j$ and $\tth_{ij}(X_{(ij)})\neq 0$. 
\item[(b)]
If $M$ is crystal module with $\rkv(M)=(1-c_{ij})\alp_i+\alp_j$ we
have $\tth_{ij}(M)=0$.
\end{itemize}
\end{Lem}
This leads us to  define in $\cM(\Pi)$ the ideal
$\cI$, which is generated by the homogeneous elements $\tth_{ij}$ for
$i,j\in I$ with $i\neq j$. 
We set moreover
\[
\bcM(\Pi)=\cM(\Pi)/\cI\ \text{and} \bar{f}:=f+\cI\quad (f\in\cM(\Pi)).
\]
Thus, by Proposition~\ref{prp:pbw1}, the morphism $\iota^*_\Ome$ induces a 
surjective algebra homomorphism 
$\bar{\iota}^*_\Ome\df\bcM(\Pi)\ra\cM(H)$. 
On the other hand, we can define for each $\br\in\NN^I$ the space of functions 
with non maximal support
\[
\cS_\br :=\{f\in\cM_\br(\Pi)\mid \dim\supp(f)<\dim H(\br)\}\text{ and }
\cS:=\oplus_{\br\in\NN^I}\cS_\br.
\]
Recall that $\dim\Pi(\br)=\dim H(\br)$. Proposition~\ref{prp:cryst} and
Lemma~\ref{lem:noSerre} imply at least that $\tth_{ij}\in\cS$. In view
of Lemma~\ref{lem:semican} and Proposition~\ref{Prp:IrrPi} it is easy
to show the following result:

\begin{Prop} The following three conditions are equivalent:
\[
(1)\ \cI\subset\cS,\qquad
(2)\ \cI=\cS,\qquad
(3)\ \cS \text{ is an ideal.}
\]
In this case the surjective algebra homomorphism
\[
\eta\df U(\fn)\ra \bcM(\Pi), e_i\mapsto \tth_i+\cI
\]
would be an isomorphism, and the  
$(\eta^{-1}(\bar{f}_Z))_{\cB}$ would form a basis of 
$U(\fn)$ which is independent of the possible choices for the $(f_Z)_{Z\in\cB}$.
\end{Prop}

Thus we call the equivalent conditions of the above proposition our
\emph{Support conjecture}.

\begin{Rem}
Our semicanonical basis would yield, similarly
to~\cite[Sec.~3]{L2}, in a natural way a basis for each integrable highest
weight representation $L(\lam)$ of $\fg(C)$, if the support conjecture is true.
See~\cite[Se.~7.3]{qrel4} for more details.
\end{Rem}
\noindent
\textbf{Acknowledgments.} I wish to thank Bernard Leclerc and Jan Schröer for 
their very precious collaboration and friendship over the years.
I also should like to acknowledge hospitality from the MPIM-Bonn, 
where this manuscript was prepared.

\end{document}